\documentclass[12pt]{article}
\usepackage[square,sort,comma,numbers]{natbib}
\usepackage{array}
\usepackage{subfigure}
\usepackage{fullpage}

\usepackage{url}
\usepackage{amsmath, amssymb, color}
\usepackage{graphicx}

\newtheorem{lemma}{Lemma}
\newtheorem{proposition}{Proposition}
\newtheorem{theorem}{Theorem}

\newcommand{\cF}{\cal F}
\newcommand{\cI}{\cal I}

\title{On Solving Convex Optimization Problems with
Linear Ascending Constraints} \author{Zizhuo Wang\thanks{Department
of Industrial and Systems Engineering, University of Minnesota,
Minneapolis, 55414. {\tt Email:zwang@umn.edu}.}}
\begin{document}
\maketitle

\abstract{In this paper, we propose two algorithms for solving
convex optimization problems with linear ascending constraints. When
the objective function is separable, we propose a dual method which
terminates in a finite number of iterations.  In particular, the
worst case complexity of our dual method improves over the
best-known result for this problem in Padakandla and Sundaresan
\cite{padakandla}. We then propose a gradient projection method to
solve a more general class of problems in which the objective
function is not necessarily separable. Numerical experiments show
that both our algorithms
work well in test problems.}\\


\section{Introduction}
\label{section:introduction} In this paper, we consider the
following optimization problem:
\begin{eqnarray}
{\bf \mbox{\bf{(P1)}}} \quad\quad\mbox{minimize}_{\vec{y}}& F(\vec{y}) = f(y_1,...,y_n) & \label{objectivemain}\\
\mbox{subject to}        & \sum_{i=1}^ky_i \le \sum_{i=1}^k \alpha_i, &\forall k = 1,...,n-1\label{firstconstraintmain}\\
& \sum_{i=1}^n y_i = (\le) \sum_{i=1}^n
\alpha_i\label{secondconstraintmain}\\
 & 0\le y_i\le \beta_i, &
\forall i = 1,...,n, \label{nonnegativeconstraintmain}
\end{eqnarray}
where $F(\cdot)$ is strictly convex in $\vec{y} = (y_1,...,y_n)$,\
and $0\le \alpha_i < +\infty$, $0\le \beta_i\le +\infty$, for $i =
1,...,n$. We make the following contributions in this paper.
\begin{enumerate}
\item We develop a dual method to solve a special case of
({\bf P1}) with separable objective functions and
(\ref{secondconstraintmain}) being an inequality constraint. Our
dual method stops in a finite number of iterations and improves the
computational complexity over the algorithm in \cite{padakandla}.
\item Using the dual method as a subroutine, we propose a gradient
projection method to solve ({\bf P1}). Our proposed method takes
advantages of the structure of the constraints so that each
projection step can be completed efficiently. The gradient
projection method also allows non-separable objective functions and
equality constraint in (\ref{secondconstraintmain}).
\item We perform numerical experiments on several test problems. The
results show that our proposed algorithms outperform the algorithm
in \cite{padakandla} as well as the standard interior point method
in most test problems.
\end{enumerate}
%

\subsection{An Alternative Form}
\label{subsec:equivalent_forms} We first point out an alternative
form of ({\bf P1}) which is sometimes used in the literature:
\begin{eqnarray}\label{equivalent}
{\bf \mbox{\bf{(P2)}}} \quad\quad\mbox{minimize}_{\vec{y}} & G(\vec{y}) = g(y_1,...,y_n) \nonumber\\
\mbox{subject to} & \sum_{i=1}^k y_i \ge \sum_{i=1}^k\alpha_i, &
\forall k =
1,...,n-1\nonumber\\
& \sum_{i=1}^n y_i = (\ge) \sum_{i=1}^n \alpha_i\nonumber \\
& 0\le y_i\le \beta_i &\forall  i = 1,...,n,
\end{eqnarray}
where $G(\vec{y})$ is strictly convex in $\vec{y}$. To translate
(\ref{equivalent}) into ({\bf P1}), we define $z_i = \beta_i - y_i$,
\footnote{It is without loss of generality to assume $\beta_i$'s are
finite, since the objective is strictly convex in $\vec{y}$,
therefore, in order to be an optimal solution, $y_i$ must be bounded
from above. The same argument applies to the other direction of
transformation.} and replace $y_i$ by $z_i$, then the optimization
problem becomes:
\begin{eqnarray*}\label{transformation}
\mbox{minimize}_{\vec{z}} & F(\vec{z}) =
G(\beta_1-z_1,...,\beta_n-z_n)\nonumber\\
\mbox{subject to} & \sum_{i=1}^k z_i \le
\sum_{i=1}^k(\beta_i-\alpha_i),& \forall k = 1,...,n-1\nonumber\\
& \sum_{i=1}^n z_i =(\le) \sum_{i=1}^n (\beta_i-\alpha_i)\nonumber\\
& 0\le z_i\le \beta_i, & \forall i=1,...,n,
\end{eqnarray*}
which is exactly of form ({\bf P1}).\footnote{It is without loss of
generality to assume that $\gamma_k =
\sum_{i=1}^k(\beta_i-\alpha_i)$ is increasing in $k$. Otherwise, we
can iteratively redefine $\gamma_k = \min_{k\le l \le n} \gamma_l$
from $n$ to $1$, and the resulting problem will be equivalent and
having the property that the right hand side of the inequality
constraints is increasing.}

\subsection{Applications}
\label{subsec:application} The formulation ({\bf P1}) arises in many
applications. One example which is a problem of smoothing is
discussed in Bellman and Dreyfus \cite{bellman}. Another one that
arises in a special case of network flow problems is studied in
Dantzig \cite{dantzig} and Veinott \cite{veinott}. Both these two
examples have the form of ({\bf P2}) with $G(\vec{y}) =
\sum_i\theta_iy_i^p$, which was also studied by Morten
\cite{morten}. Other problems of such a form arise frequently in
communication networks and are discussed in e.g., Padakandla and
Sundaresan \cite{padakandla,padakandla2} and Viswanath and
Anantharam \cite{visawanath}. In addition to the above applications,
we present another motivation of this model in operations
management.

{\bf Inventory problem with downward substitution.} A firm sells a
product with $n$ different grades, with $1$ the highest and $n$ the
lowest. The firm has $\alpha_i$ grade $i$ products on hand and the
demand of grade $i$ product is a random variable $D_i$. Any product
of grade $i$ can be used to satisfy the demand of product of grade
$i$ or lower ($j\ge i$). Before the demand realizes, the firm has to
make an inventory decision $y_i$ of how many grade $i$ products to
put into stock. Once this is done, the products are no longer
substitutable (for example, the firm has to package these products
during this process, products of different grades need different
packages and will not be distinguishable after packaging). For each
grade $i$ product, there is a unit overage cost $o_i$ if $D_i$ turns
out to be less than $y_i$ and a unit underage cost $u_i$ if $D_i$
turns out to be greater than $y_i$. The objective is to minimize the
expected total cost. The problem can be written as (we use ${\mathbb
E}(\cdot)$ to denote the expectation operator):
\begin{eqnarray}\label{inventory}
\mbox{minimize}_{\vec{y}} & \sum_{i=1}^n \left(u_i{\mathbb
E}(D_i-y_i)^+ + o_i
{\mathbb E}(y_i-D_i)^+\right)\nonumber\\
\mbox{subject to} & \sum_{i=1}^k y_i \le \sum_{i=1}^k \alpha_i, &
\forall k
= 1,...,n\nonumber\\
& y_i\ge 0, & \forall i = 1,...,n.
\end{eqnarray}
Note that (\ref{inventory}) is in the form of ({\bf P1}). If the
cumulative distribution function of each $D_i$ is continuous and
strictly increasing on the positive domain, the objective function
of (\ref{inventory}) is also strictly convex and differentiable. As
pointed out in the literature, such inventory problems with downward
substitution possibility occur in several practical settings such as
in semiconductor industry where higher quality chips can be used to
substitute lower quality ones, see \cite{hsu}, \cite{bassok},
\cite{rao} and \cite{wagner}. For example, Bassok et al.
\cite{bassok} study the optimal single-period inventory decision in
a production system with downward substitution. Hsu and Bassok
\cite{hsu} further incorporate random yield in such systems and
decide the optimal initial production quantity. Rao et al.
\cite{rao} further extend it to consider the problem with setup
costs. In most of these literature, the substitution is modeled to
occur after the demand is realized and is used as a recourse action,
therefore, the optimal substitution decision can be solved from a
linear program. However, one can easily envision that in some
practical problems, substitution decision has to be made prior to
the demand realization. For example, the firm may need to package
different products before delivering to the markets. In such cases,
(\ref{inventory}) will be a more appropriate model.

In practice, one has strong incentive to solve (\ref{inventory})
faster. For example, the firms may also need to decide the upfront
production quantities $\alpha_i$'s for each grade with production
cost $c_i\alpha_i$. The actual yield of grade $i$ product is
$\alpha_iU_i$ where $U_i$'s are some known random yield
distributions. In this case, the firm's problem is the following
two-stage stochastic programming problem:
\begin{eqnarray}\label{twostagesp}
\mbox{minimize}_{\vec{\alpha}} \sum_{i=1}^n c_i\alpha_i + & {\mathbb
E}_U\{\mbox{minimize}_{\vec{y}}& \sum_{i=1}^n \left(u_i{\mathbb
E}(D_i-y_i)^+ + o_i
{\mathbb E}(y_i-D_i)^+\right)\}\nonumber\\
&\mbox{subject to} & \sum_{i=1}^k y_i \le \sum_{i=1}^k \alpha_iU_i,
\quad\quad \forall k = 1,...,n\nonumber\\
&& y_i\ge 0, \quad\quad \forall i = 1,...,n.
\end{eqnarray}
This is a similar problem as introduced in \cite{hsu} except the
substitution decision has to be made before demand realization. A
natural approach to solve (\ref{twostagesp}) is to use the
stochastic gradient method \cite{shapiro}. However, this requires
one to evaluate the inside problem repeatedly. Therefore, improving
the efficiency of solving (\ref{inventory}) could be of strong
interest.

\subsection{Literature Review}
The main related literature to this paper is \cite{padakandla}. In
\cite{padakandla}, the authors propose a dual method for solving
({\bf P2}) with separable objective functions. We call the algorithm
in \cite{padakandla} the ``P-S algorithm'' in the rest of the
discussions. The P-S algorithm is currently the state-of-the-art
algorithm for solving this problem. It finishes in $O(n)$ outer
iterations. In each iteration, it solves up to $n$ nonlinear
equations, and sets at least one primal variable based on the
solutions to the equations. The efficiency of the P-S algorithm
depends on how fast one can solve those equations. When the
equations have closed form solutions, the P-S algorithm performs
very well, otherwise, it may not. In this paper, we propose a dual
algorithm which does not attempt to set primal variables in each
iteration. Instead, we set one dual variable in each iteration and
maintain the optimality conditions for the variables that have been
set. Our dual algorithm also finishes in $O(n)$ outer iterations and
in each iteration, we solve no more than one equation. We show that
the equations we solve are simply the equations in the P-S algorithm
with lower bound on each term. When the equations in the P-S
algorithm do not have a closed form solution, solving both equations
usually have the same complexity. In those cases, our dual algorithm
reduces the computational complexity of the P-S algorithm by an
order of $n$.

In addition to the dual method, we propose a gradient projection
method to solve the more general problem ({\bf P1}) allowing
non-separable objective functions. Gradient projection methods are
widely used to solve a variety of convex optimization problems. We
refer the readers to \cite{bertsekas} for a thorough discussion of
this method. In particular, the key element in the gradient
projection method is the design of the projection step. In this
paper, we propose an efficient projection step under linear
ascending constraints which leads to an efficient implementation of
the gradient projection algorithm for the problem.

Another popular method that solves nonlinear convex optimization is
the interior point method. However, we focus on the first order
method in this paper because of its low memory requirement and thus
the ability to solve large problems. Performance comparisons between
our proposed algorithms and the interior point algorithm
(implemented by CVX) are shown in the numerical tests and the
results indicate that our algorithms are usually much more
efficient.

We note that there is abundant literature on solving a special case
of ({\bf P1}) when there is only one equality/inequality constraint
(usually called the simplex constraint or $l_1$ constraint). We
refer the readers to \cite{patriksson} for a survey on this problem.
Although the dual method is widely used in those studies, the detail
of our algorithm differs significantly because of the special
structure of this problem.

\subsection{Structure of the paper}

In Section \ref{section:dual method}, we develop a dual method to
solve a special case of ({\bf P1}) with separable objective
functions and (\ref{secondconstraintmain}) being an inequality
constraint. In Section \ref{section:gradient method}, we further
propose a gradient projection method to solve the general problem
({\bf P1}). Numerical tests are shown in Section
\ref{section:numerical_experiments} to examine the performances of
our algorithms. Section \ref{section:conclusion} concludes this
paper.

\section{A Dual Method}
\label{section:dual method}

In this section, we study a special case of ({\bf P1}) in which the
objective function $F$ is separable, i.e., $F(\vec{y}) =
\sum_{i=1}^n f_i(y_i)$ and (\ref{secondconstraintmain}) is an
inequality constraint. There are two reasons why we consider
separable objectives. First, in most of the applications mentioned
in Section \ref{subsec:application}, the objective functions are
indeed separable. Second, the study of separable objective functions
will lay the foundation for the analysis of the gradient projection
method in Section \ref{section:gradient method} which can solve more
general problems.

In the following, we develop a dual method to solve the following
problem:
\begin{eqnarray}\label{model_separable}
{\bf \mbox{\bf{(P3)}}} \quad\quad\mbox{minimize}_{\vec{y}}&
\sum_{i=1}^nf_i(y_i) &
\label{objective}\\
\mbox{subject to}        & \sum_{i=1}^ky_i \le \sum_{i=1}^k \alpha_i, &\forall k = 1,...,n\label{sumconstraint}\\
& 0\le y_i\le \beta_i, & \forall i = 1,...,n.
\label{nonnegativeconstraint}
\end{eqnarray}
Here we assume that $f_i(y_i)$'s are continuously differentiable,
strictly convex in $y_i$ with derivative $g_i(y_i) = f_i'(y_i)$.
\footnote{Our algorithm works in a similar manner even if $f$ is not
differentiable but convex. The discussions will involve subgradient
of $f$ in that case. We make this assumption simply for the
convenience of discussion.} Under these assumptions, $g_i(\cdot)$ is
strictly increasing and we denote $g_i(0) = l_i$ and $g_i(\beta_i) =
h_i$. We define $\bar{y}_i = arg\min_{0\le y\le \beta_i} f_i(y)$,
that is, $\bar{y}_i$'s are the optimal solution to ({\bf P3})
without constraint (\ref{sumconstraint}). Under the above
assumptions, it is easy to see that $\bar{y}_i$ exists and is
unique.

We first write down the KKT conditions of
(\ref{objective})-(\ref{nonnegativeconstraint}). We associate a dual
variable $\lambda_k$ to each constraint (\ref{sumconstraint}), a
dual variable $\delta_i$ to each upper bound constraint, and a dual
variable $\eta_i$ to each nonnegative constraint
(\ref{nonnegativeconstraint}). The Lagrangian of
(\ref{objective})-(\ref{nonnegativeconstraint}) can then be written
as
\begin{eqnarray*}
\sum_{i=1}^nf_i(y_i) +\sum_{k=1}^n \lambda_k\left(\sum_{i=1}^k y_i -
\sum_{i=1}^k \alpha_i\right) - \sum_{i=1}^n \eta_i y_i -
\sum_{i=1}^n \delta_i(\beta_i-y_i).
\end{eqnarray*}
And the KKT conditions are
\begin{eqnarray}
g_i(y_i) = -\sum_{k=i}^n \lambda_k+\eta_i - \delta_i, \quad\quad\forall i=1,...,n, \label{kkt1}\\
y_i\cdot \eta_i =0, y_i\ge 0, \eta_i\ge 0, \quad\quad\forall i=1,...,n, \label{kkt2} \\
(\beta_i -y_i) \cdot\delta_i = 0, y_i\le \beta_i, \delta_i\ge 0,
\quad\quad\forall i=1,...,n, \label{kkt3}\\ \sum_{i=1}^k y_i \le
\sum_{i=1}^k \alpha_i, \quad\quad\forall k=1,...,n,
\label{kkt4}\\
\lambda_k \cdot \left(\sum_{i=1}^ky_i -\sum_{i=1}^k\alpha_i\right) =
0, \lambda_k \ge 0, \quad\quad\forall k=1,...,n. \label{kkt5}
\end{eqnarray}
Define $\phi_i(x) = \max\{l_i, \min\{x, h_i\}\}$, i.e., $\phi_i(x)$
projects $x$ to the interval $[l_i, h_i]$. We also define $H_i(x) =
g_i^{-1}(\phi_i(x))$ where $g_i^{-1}(\cdot)$ is the inverse function
of $g_i(\cdot)$. By the assumptions on $g_i(\cdot)$, $l_i$ and
$h_i$, we have $0\le H_i(x) \le\beta_i$. From
(\ref{kkt1})-(\ref{kkt3}), one can observe that $g_i(y_i)$ should
equal to the projection of $-\sum_{k=i}^n \lambda_k$ onto the
interval $[l_i,h_i]$. More precisely, conditions
(\ref{kkt1})-(\ref{kkt3}) can be equivalently written as
\begin{eqnarray*}
y_i = H_i\left(-\sum_{k=i}^n\lambda_k\right)
\end{eqnarray*}
with
\begin{eqnarray}
\eta_i = \left(\phi_i\left(-\sum_{k=i}^n\lambda_k\right) +
\sum_{k=i}^n\lambda_k\right)^+\mbox{  and } \delta_i =
\left(-\phi_i\left(-\sum_{k=i}^n\lambda_k\right)
-\sum_{k=i}^n\lambda_k\right)^+, \forall i \label{kkt6}
\end{eqnarray}
where  $x^+ = \max\{x,0\}$. Since {\bf{(P3)}} is linearly
constrained and is convex, the KKT conditions are necessary and
sufficient, and solving them yields the solution of {\bf{(P3)}}. In
the following, we propose an efficient dual method to solve the KKT
conditions. The idea of this dual method is to assign values to the
dual variables $\lambda$'s such that the optimality conditions
(\ref{kkt4})-(\ref{kkt6}) hold. We state our algorithm as follows:

\hrulefill\\
{\noindent\bf Algorithm 1}\\

{\noindent\bf Step 0: Initialization.} Let $d_k =
\sum_{i=1}^k\bar{y}_i - \sum_{i=1}^k \alpha_i$, $k = 1,2,...,n$.
Define
\begin{eqnarray*}
w_0 & = & 0,\\
w_1 & =&  \min\{k: d_k\ge 0\},\\
w_{j+1} &= &\min\{k > w_j: d_k\ge d_{w_{j}}\},\mbox{ } \forall j.
\end{eqnarray*}
Here we define $\min \emptyset = \infty$. If $w_1 = \infty$, then
setting $y_i = \bar{y}_i = arg\min_{0\le y\le \beta_i} f_i(y)$ and
$\lambda_i = 0$ for all $i$ will satisfy the KKT conditions and thus
is optimal. Otherwise, let $L = \max\{j\ge 1: w_j<\infty\}$. Define
$S= \{w_1,w_2,...,w_L\}$. Let $\lambda_i =
0$, $\eta_i =0$ for all $i$ and let $j = L$.\\

{\noindent\bf Step 1: Main Loop (Outer Loop)}. \\
{\bf WHILE} $j > 0$ 
\begin{itemize}
\item Case 1: If
\begin{eqnarray}\label{iter1case1}
\sum_{s=w_{j-1}+1}^{w_j}\left(H_s\left(- \sum_{t=j+1}^L
\lambda_{w_t}\right) - \alpha_s\right) \ge 0
\end{eqnarray}
then choose $\xi\ge 0$ such that
\begin{eqnarray}\label{iter1case1'}
\sum_{s=w_{j-1}+1}^{w_j}\left(H_s\left(- \sum_{t=j+1}^L
\lambda_{w_t} - \xi\right) - \alpha_s\right) = 0
\end{eqnarray}
and set $\lambda_{w_j} = \xi$, $j \leftarrow j - 1$.
\item Case 2: If
\begin{eqnarray}\label{iter1case2}
\sum_{s=w_{j-1}+1}^{w_j}\left(H_s\left(- \sum_{t=j+1}^L
\lambda_{w_t}\right) - \alpha_s\right) < 0
\end{eqnarray}
then use binary search to find
\begin{eqnarray}\label{r}
r^* = \min_r\left\{j+1\le r \le L:
\sum_{s={w_{j-1}+1}}^{w_r}\left(H_s\left(- \sum_{t=r+1}^L
\lambda_{w_t}\right)-\alpha_s\right) \ge 0\right\}.
\end{eqnarray}
If such $r^*$ does not exist, then set all $\lambda_{w_r} = 0$, for
$r = j,...,L$. Set $j \leftarrow j - 1$. Otherwise, choose $\xi\ge
0$ such that
\begin{eqnarray}\label{iter2case1'}
\sum_{s={w_{j-1}+1}}^{w_{r^*}}\left(H_s\left( -\sum_{t=r^*+1}^L
\lambda_{w_t} -\xi\right)-\alpha_s\right) = 0
\end{eqnarray}
and set $\lambda_{w_{r^*}} = \xi$ and $\lambda_{w_r} = 0$, for $j\le
r < r^* $. Set $j \leftarrow j - 1$.
\end{itemize}
{\bf END WHILE}\\

{\noindent\bf Step 2: Output} Set for $i=1,2,...,n$
\begin{eqnarray*} y_i =
H_i\left(-\sum_{k=i}^n\lambda_k\right),
\end{eqnarray*}
\begin{eqnarray*}
\eta_i = \left(\phi_i\left(-\sum_{k=i}^n\lambda_k\right) +
\sum_{k=i}^n\lambda_k\right)^+ \quad\mbox{  and  }\quad \delta_i =
\left(-
\phi_i\left(-\sum_{k=i}^n\lambda_k\right)-\sum_{k=i}^n\lambda_k\right)^+.
\end{eqnarray*}
\hrulefill\\

First, we argue that those $\xi$'s defined in (\ref{iter1case1'})
and (\ref{iter2case1'}) exist. This can be verified by observing
that when $\xi = 0$, the left hand sides of (\ref{iter1case1'}) and
(\ref{iter2case1'}) are both nonnegative and as $\xi \rightarrow
\infty$, both of them will be less than or equal to zero (since
$\alpha$'s are all nonnegative). Also, by our assumption, the left
hand sides of (\ref{iter1case1'}) and (\ref{iter2case1'}) are both
continuous. Therefore, by the intermediate value theorem, such
$\xi$'s must exist. We now state the main result of this section.

\begin{theorem}
\label{theorem:dual} Algorithm 1 terminates within $L\le n$ outer
iterations and the output solves {\bf{(P3)}}.
\end{theorem}

In the following, we prove Theorem \ref{theorem:dual}. Let
$(\{y_i^*\}_{i=1}^n,\{\lambda_i^*\}_{i=1}^n,\{\eta_i^*\}_{i=1}^n,
\{\delta_i^*\}_{i=1}^n)$ be any solution to the KKT conditions
(\ref{kkt4})-(\ref{kkt6}) and thus an optimal solution to
{\bf{(P3)}}. First, it is easy to see that $y_i^*\le \bar{y}_i$ for
all $i$, otherwise replacing $y_i^*$ with $\bar{y}_i$ will strictly
improve the objective value while still satisfying the constraints,
which contradicts with the optimality of $y_i^*$'s. Next we claim
that for any $k\in \{1,...,n\}\setminus S$, we must have
$\lambda_k^* = 0$. This is because for  such $k$ with
$w_{l-1}<k<w_l$, we have
\begin{eqnarray*}
\sum_{i=1}^{k} (y_i^*-\alpha_i) \le \sum_{i=w_{l-1}+1}^k
(y_i^*-\alpha_i) \le \sum_{i=w_{l-1}+1}^{k}(\bar{y}_i-\alpha_i) < 0,
\end{eqnarray*}
where the last inequality is because of the definition of $w_l$. By
the complementarity condition (\ref{kkt5}), $\lambda_k^* = 0$ for
$w_{l-1}<k<w_l$.

Note that in the KKT conditions, given $\lambda_k$'s, the $y$'s,
$\eta$'s and $\delta$'s are uniquely determined and that changing
$\lambda_k$ only affects $y_i$'s, $\eta_i$'s and $\delta_i$'s for $i
\le k$. In iteration $j$ of Algorithm 1, we assign $\lambda_{w_j}$
and may modify all $\lambda_{w_k}$'s for $k\ge j+1$. We now state
the following property of Algorithm 1 which immediately implies
Theorem \ref{theorem:dual}.

\begin{proposition}\label{claim:1}
When Algorithm 1 finishes loop $j$ ($ j = L, L-1,..., 1$), the
current $\lambda_i$'s together with

\begin{eqnarray}\label{currenty}
y_i = H_i\left(-\sum_{k=i}^n\lambda_k\right),
\end{eqnarray}
\begin{eqnarray}\label{currenteta}
\eta_i = \left(\phi_i\left(-\sum_{k=i}^n\lambda_k\right) +
\sum_{k=i}^n\lambda_k\right)^+ \quad\mbox{  and  }\quad \delta_i =
\left(-
\phi_i\left(-\sum_{k=i}^n\lambda_k\right)-\sum_{k=i}^n\lambda_k\right)^+
\end{eqnarray}

satisfy the following conditions:
\begin{eqnarray}\label{tempkkt1}
y_i\cdot \eta_i = 0, \eta_i \ge 0 , y_i\ge 0 ,\quad\quad\forall i\ge
w_{j-1}+1,
\end{eqnarray}
\begin{eqnarray}\label{tempkkt1.5}
(\beta_i-y_i)\cdot \delta_i = 0, y_i\le \beta_i, \delta_i\ge 0,
\quad\quad\forall i \ge w_{j-1}+1,
\end{eqnarray}
\begin{eqnarray}\label{tempkkt3}
\sum_{s=w_{j-1}+1}^k y_s\le \sum_{s=w_{j-1}+1}^k
\alpha_s,\quad\forall k \ge w_{j-1} + 1,
\end{eqnarray}
\begin{eqnarray}\label{tempkkt2}
\lambda_k\cdot\left(\sum_{s=w_{j-1}+1}^k y_s - \sum_{s=w_{j-1}+1}^k
\alpha_s\right) = 0, \quad\forall k\ge w_{j-1} + 1.
\end{eqnarray}
\end{proposition}

Before we prove Proposition \ref{claim:1}, we introduce a lemma that
will be used repeatedly in the proof.
\begin{lemma}\label{lemma2}
For all $i=1,...,n$, if $y_i$ is defined in (\ref{currenty}), then
$y_i\le \bar{y}_i$.
\end{lemma}
The lemma follows immediately from the assumption that
$g_i(\cdot)$'s are strictly increasing and that $\bar{y}_i =
H_i(0)$.\\

{\bf\noindent Proof of Proposition \ref{claim:1}:} First, note that
condition (\ref{tempkkt1}) and (\ref{tempkkt1.5}) are always
satisfied because of the definitions in (\ref{currenty}) and
(\ref{currenteta}). Therefore, it suffices to show that conditions
(\ref{tempkkt3}) and (\ref{tempkkt2}) hold for $j = L, L-1,...,1$.
We use backward induction to prove this. First we show that for
$j=L$, (\ref{tempkkt3}) and (\ref{tempkkt2}) hold for all $k \ge
w_{L-1}+1$.

First we show that (\ref{tempkkt3}) holds. When Algorithm 1 finishes
the outer loop when $j = L$, for any $w_{L-1}+1\le s' <w_L$, we have
\begin{eqnarray}
\sum_{s=w_{L-1}+1}^{s'}(y_s-\alpha_s)\le \sum_{s=w_{L-1}+1}^{s'}
(\bar{y}_s -\alpha_s) \le 0,
\end{eqnarray}
where the first inequality is due to Lemma \ref{lemma2} and the
second inequality is due to the definition of $w_L$. On the other
hand, for $s'\ge w_L$, we have
\begin{eqnarray*}
\sum_{s=w_{L-1}+1}^{s'} (y_s-\alpha_s) \le \sum_{s=w_{L-1}+1}^{w_L}
(y_s-\alpha_s) + \sum_{s=w_L + 1}^{s'} (\bar{y}_s-\alpha_s) \le  0,
\end{eqnarray*}
where the first inequality is because of Lemma \ref{lemma2} and the
second one is because of step (\ref{iter1case1'}) and the definition
of $w_L$. Therefore (\ref{tempkkt3}) holds when $j = L$.

To show that (\ref{tempkkt2}) holds for $j = L$, note that among all
the $\lambda_k$'s with $k\ge w_{L-1}+1$, the only possible non-zero
one is $\lambda_{w_L}$. If Case 1 of the algorithm happens in this
loop, then
\begin{eqnarray*}
\sum_{s=w_{L-1}+1}^{w_L} (y_s-\alpha_s) = 0.
\end{eqnarray*}
Otherwise, $\lambda_{w_L} = 0$. Therefore, (\ref{tempkkt2}) holds
for $j = L$.

Now we assume that (\ref{tempkkt3}) - (\ref{tempkkt2}) hold after
the algorithm completes the outer loop for $j = \bar{j} + 1$. Now we
consider the situation when it finishes the outer loop for
$j=\bar{j}$. We consider two cases:

\begin{itemize}
\item Case 1: (\ref{iter1case1}) holds in the current loop ($j = \bar{j}$). In this case, we
have
\begin{eqnarray}\label{16}
\sum_{s=w_{\bar{j}-1}+1}^{w_{\bar{j}}} y_s  =
\sum_{s=w_{\bar{j}-1}+1}^{w_{\bar{j}}} \alpha_s.
\end{eqnarray}
And the $y_s$'s for $s > w_{\bar{j}}$ does not change from the
previous loop. Therefore, for any $k = w_j$ $(j\ge \bar{j})$, we
have
\begin{eqnarray*}
\sum_{s = w_{\bar{j}-1}+1}^k\ y_s \le \sum_{s=w_{\bar{j}-1}+1}^k
\alpha_s.
\end{eqnarray*}
And for $w_j < k < w_{j+1}$ ($j\ge \bar{j} - 1$),
\begin{eqnarray*}
\sum_{s=w_{\bar{j}-1}+1}^k (y_s-\alpha_s) & \le &
\sum_{s=w_{\bar{j}-1}+1}^{w_j} (y_s-\alpha_s) + \sum_{s=w_j+1}^k
(\bar{y}_s-\alpha_s) \le \sum_{s = w_{\bar{j}-1}+1}^{w_j}
(y_s-\alpha_s)\le 0,
\end{eqnarray*}
where the first inequality is because of Lemma \ref{lemma2} and the
second inequality is because of the definition of $w_j$'s.
Therefore, (\ref{tempkkt3}) holds for all $k\ge w_{\bar{j}-1}+1$.
For (\ref{tempkkt2}), we only need to study $k = w_j$ since all
other $\lambda_k$'s are $0$. And it holds because of (\ref{16}) and
the induction assumption. Therefore, (\ref{tempkkt3}) -
(\ref{tempkkt2}) hold for $\bar{j}$ in this case.

\item Case 2: (\ref{iter1case2}) holds in the current loop. Then
there are two further cases:
\begin{itemize}
\item a): $r^*$ does not exist. In
this case, by the definition of Algorithm 1, all $\lambda_k$'s are
zero at the end of this iteration and $\sum_{s = w_{\bar{j}-1}+1}^k\
y_s < \sum_{s=w_{\bar{j}-1}+1}^k \alpha_s$ for all $k = w_j$ ($j\ge
\bar{j}$). By the same argument as in case 1, we know that $\sum_{s
= w_{\bar{j}-1}+1}^k\ y_s < \sum_{s=w_{\bar{j}-1}+1}^k \alpha_s$ for
all $k \ge w_{\bar{j}-1}+1$. Therefore, (\ref{tempkkt3}) -
(\ref{tempkkt2}) hold for $\bar{j}$ in this case.

\item b): $r^*$ exists. Denote the $\lambda$'s and $y$'s after the previous
loop by $\tilde{\lambda}$ and $\tilde{y}$. It is easy to see that in
this case, at the end of the current iteration, $\lambda_i \le
\tilde{\lambda}_i$ and $y_i \ge \tilde{y}_i$ for all $i$. We first
have the following lemma whose proof is relegated to Appendix
\ref{appendix:proofoflemma1}:
\begin{lemma}\label{lemma1}
$\tilde{\lambda}_{w_{r^*}} > 0$.
\end{lemma}

With Lemma \ref{lemma1}, we show that (\ref{tempkkt3}) -
(\ref{tempkkt2}) hold. We first consider (\ref{tempkkt3}). By
(\ref{iter2case1'}), we have
\begin{eqnarray}\label{18}
\sum_{s = w_{\bar{j}-1}+1}^{w_{r^*}} (y_s-\alpha_s) = 0.
\end{eqnarray}
Therefore for $k > w_{r^*}$, we know that
\begin{eqnarray*}
\sum_{s=w_{\bar{j}-1}+1}^k (y_s-\alpha_s) = \sum_{s=w_{r^*}+1}^k
(y_s-\alpha_s) = \sum_{s=w_{r^*}+1}^k (\tilde{y}_s-\alpha_s) =
\sum_{s=w_{\bar{j}}+1}^k (\tilde{y}_s-\alpha_s) \le 0,
\end{eqnarray*}
where the second equality is because the value of $y_s$ does not
change for $s\ge w_{r^*}+1$. The last equality is because of the
induction assumption that $\tilde{\lambda}_{w_{r^*}}
\cdot\sum_{s=w_{\bar{j}}+1}^{w_{r^*}}(\tilde{y}_s-\alpha_s) = 0$ and
Lemma \ref{lemma1}. Therefore, for all $w_{j}\le k < w_{j+1}$,
$\bar{j}-1\le j < r^*$, we have
\begin{eqnarray*}
\sum_{s=w_{\bar{j}-1}+1}^k(y_s-\alpha_s)\le
\sum_{s=w_{\bar{j}-1}+1}^{w_j} (y_s-\alpha_s)
 \le -\sum_{s=w_j+1}^{w_{r^*}} (\tilde{y}_s-\alpha_s) = \sum_{s=
w_{\bar{j}}+ 1}^{w_j} (\tilde{y}_s-\alpha_s) \le 0,
\end{eqnarray*}
where the first inequality is because of the definition of $w_j$,
the second equality is because of (\ref{18}) and $y_i\ge
\tilde{y}_i$ and the last equality is because of the induction
assumption and Lemma \ref{lemma1}. Therefore, (\ref{tempkkt3}) holds
in this case.

Lastly, we show that (\ref{tempkkt2}) also holds. It suffices to
show that for each $r > r^*$ such that $\lambda_r > 0$,
\begin{eqnarray*}
\sum_{s= w_{\bar{j}-1}+1}^{w_r}(y_s-\alpha_s) = 0.
\end{eqnarray*}
This is equivalent as showing that for each $r > r^*$ such that
$\lambda_r > 0$,
\begin{eqnarray*}
\sum_{s=w_r^* + 1}^{w_r} (y_s-\alpha_s) = 0.
\end{eqnarray*}
Note that
\begin{eqnarray*}
\sum_{s=w_{r^*} + 1}^{w_r} (y_s-\alpha_s) = \sum_{s=w_{r^*} +
1}^{w_r} (\tilde{y}_s-\alpha_s) =
\sum_{s=w_{\bar{j}}+1}^{w_r}(\tilde{y}_s-\alpha_s) -
\sum_{s=w_{\bar{j}}+1}^{w_{r^*}}(\tilde{y}_s-\alpha_s)
\end{eqnarray*}
By induction assumption and Lemma \ref{lemma1},
\begin{eqnarray*}
\sum_{s=w_{\bar{j}}+1}^{w_r}(\tilde{y}_s-\alpha_s) =
\sum_{s=w_{\bar{j}}+1}^{w_{r^*}}(\tilde{y}_s-\alpha_s) = 0.
\end{eqnarray*}
Therefore (\ref{tempkkt2}) holds in this case and Proposition
\ref{claim:1} is proved. $\hfill \Box$\\
\end{itemize}
\end{itemize}

Now we make some comments on Algorithm 1.

By its definition, Algorithm 1 terminates within $L\le n$ outer
iterations. In practical problems, $L$ might be much less than $n$.
In those cases, the algorithm can output the solutions very fast.
This is a similar property as in the P-S algorithm (recall we use
P-S algorithm to refer the algorithm proposed in \cite{padakandla}).
Now we use $\cI$ to denote the complexity (number of arithmetic
operations) of solving (\ref{iter1case1'}) or (\ref{iter2case1'})
once (it is easy to see that ${\cI}\ge n$). In each iteration of
Algorithm 1, if Case 1 happens, the algorithm has to perform a sum
of no more than $n$ terms. And it has to solve (\ref{iter1case1'})
once.  Therefore, there are $O(\cI)$ arithmetic operations in this
case. If Case 2 happens, then the algorithm has similar tasks as in
Case 1, and in addition it needs to find $r^*$ defined in (\ref{r})
which takes no more than $O(n\log{n})$ iterations. Therefore, the
complexity in Case 2 is $O(\max\{n\log{n}, {\cI}\})$. Combined with
$O(n)$ outer iterations , the total arithmetic complexity of our
algorithm is $O(\max\{n^2\log{n}, n{\cI}\})$.

Now we compare the complexity result to that of the P-S algorithm.
The difference between the two algorithms is the way the variables
are assigned. In each iteration of the P-S algorithm, it solves
\begin{eqnarray}\label{refequation}
\sum_{m\in s\cap [i,l]} (g_m^{-1}(\theta)\wedge \beta_m) =
\sum_{m=i}^l \alpha_m
\end{eqnarray}
for all $l \le j$, where $s$ is the set of unassigned variables.
Then the largest solution is chosen and the corresponding primal
variable is set accordingly. Such a method avoids the needs to check
the validity of the KKT conditions that is met in previous steps as
we have to do in Step 2 in Algorithm 1, however at a cost of having
to solve $O(n)$ equations at each step rather than only one as in
Algorithm 1. Indeed, the equations (\ref{refequation}) are sometimes
easier to solve since they don't involve the lower bound as
Algorithm 1 do. If one denotes the arithmetic complexity of solving
equations in (\ref{refequation}) by $O({\cI}')$, then the total
arithmetic complexity of the P-S algorithm is
$O(n^2{\cI}')$.\footnote{Again, it is easy to see that ${\cI}'$ is
at least $O(n)$ since one has to sum $O(n)$ values in order to solve
(\ref{refequation}).} Therefore, our algorithm works better than the
P-S algorithm when solving equations in (\ref{refequation}) has
similar complexity as solving equations in (\ref{iter1case1'}) and
(\ref{iter2case1'}), but may work relatively worse if
(\ref{refequation}) can be solved explicitly (see \cite{padakandla}
for several examples). This tradeoff is demonstrated in the
numerical experiments in Section
\ref{section:numerical_experiments}.

There are two main drawbacks for Algorithm 1. First, it can only
handle separable objective functions and inequality constraint in
(\ref{secondconstraintmain}). Second, it involves many evaluations
of $g^{-1}$ and also has to solve the equations (\ref{iter1case1'})
and (\ref{iter2case1'}). These evaluations might be very expensive
in computation if $g^{-1}$ does not have a simple form. This is the
same problem as in the P-S algorithm. In particular,
\cite{padakandla} shows that the performance of the  P-S algorithm
may not be very good if closed form solutions to equations
(\ref{refequation}) do not exist. To overcome this problem, we
propose a gradient projection method in the next section. The
gradient projection method uses Algorithm 1 as a subroutine,
however, in each iteration, $g(\cdot)$ is simply a linear function.
Moreover, the gradient projection method can handle non-separable
objective functions as well as equality constraints in
(\ref{secondconstraintmain}). The tradeoff however, is that the
gradient projection method does not give an exact solution in a
finite number of iterations. However as we demonstrate in our
numerical experiments, it performs quite well in test problems.

\section{Gradient Projection Method}
\label{section:gradient method}

In this section, we propose a gradient projection method to solve
({\bf P1}). First we claim that we can assume that constraint
(\ref{secondconstraintmain}) is of the inequality form. To transform
a problem with equality constraint in (\ref{secondconstraintmain})
to an inequality one, we first note that we can without loss of
generality assume $\beta_n = \infty$. This is because one can always
add a penalty term $M(y_n-\beta_n)^+$ with sufficiently large $M$ so
that the optimal solution must satisfy $y_n\le \beta_n$ (if the
problem is feasible). Then, we can simply substitute $y_n =
\sum_{i=1}^n \alpha_i - \sum_{i=1}^{n-1} y_i$ into
(\ref{objective}). Therefore, it is sufficient to consider the
following equivalent problem:
\begin{eqnarray}\label{model3}
{\bf \mbox{\bf{(P4)}}} \quad\quad\mbox{minimize}_{\vec{y}}& F(\vec{y}) = f(y_1,...,y_n) & \nonumber\\
\mbox{subject to}        & \sum_{i=1}^ky_i \le \sum_{i=1}^k \alpha_i, &\forall k = 1,...,n\nonumber\\
 & 0\le y_i\le \beta_i, & \forall i = 1,...,n.
\end{eqnarray}

In the following, we propose a gradient projection method to solve
{\bf(P4)}. Gradient projection methods are used to solve a variety
of convex optimization problems \cite{bertsekas}. It minimizes a
function $F(x)$ subject to convex constraints by generating the
sequence $\left\{x^{(k)}\right\}$ via
\begin{eqnarray*}\label{grdient_projection_iter}
x^{(k+1)} = \Pi^k\left(x^{(k)} - \mu_k\nabla^{(k)}\right),
\end{eqnarray*}
where $\nabla^{(k)}$ is the gradient of $F(x)$ at $x^{(k)}$, $\Pi(x)
= arg\min_{y}\left\{||x-y||:y\in {\cF}\right\}$ is the Euclidean
projection of $x$ onto the feasible set ${\cF}$ and $\mu_k$ is the
chosen stepsize. In the following, our discussion will mainly focus
on the projection step. The convergence of the gradient projection
method is referred to \cite{bertsekas}.

In our problem, given $z = x^{(k)} - \mu_k\nabla^{(k)}$, $x^{(k+1)}$
can be computed by solving
\begin{eqnarray}\label{projection}
\mbox{minimize}_{\vec{y}}& \sum_{i=1}^n(y_i- z_i)^2 & \nonumber\\
\mbox{subject to}        & \sum_{i=1}^ky_i \le \sum_{i=1}^k \alpha_i, &\forall k = 1,...,n\nonumber\\
 & 0\le y_i\le \beta_i, & \forall i = 1,...,n.
\end{eqnarray}
Note that (\ref{projection}) is of form {\bf (P3)} thus can be
solved by Algorithm 1. One main advantage of (\ref{projection}) is
that the objective function is quadratic. Therefore, in Algorithm 1,
$g_i(y_i) = 2(y_i-z_i)$, $l_i = -2z_i$, $h_i = 2(\beta_i - z_i)$ and
$g_i^{-1}(u_i) = \frac{u_i + 2z_i}{2}$. Therefore the equation
(\ref{iter1case1'}) and similarly (\ref{iter2case1'}) can be written
as
\begin{eqnarray}\label{equationinprojection}
\rho(\xi) = \sum_{s=w_{j-1}+1}^{w_j}
\left(\frac{\max(0,\min(2\beta_i,2z_s-\sum_{t=j+1}^L\lambda_{w_t}-\xi))}{2}-\alpha_s\right)
= 0.
\end{eqnarray}
Note that (\ref{equationinprojection}) is a decreasing piecewise
linear function with no more than $2n$ breakpoints. And those
breakpoints can be computed explicitly. Therefore, to solve
(\ref{equationinprojection}), one can first use binary search to
find out which piece of the function the solution belongs to and
then simply solve a linear equation. Therefore, the total complexity
of solving (\ref{equationinprojection}) is $O(n\log{n})$ and the
total complexity of each projection step is $O(n^2\log^2{n})$,
regardless of the form of the objective function.


\section{Numerical Experiments}
\label{section:numerical_experiments}

In this section, we perform numerical tests to examine the
performance of both our dual method and the gradient projection
method and compare them to 1) the P-S algorithm in \cite{padakandla}
and 2) CVX \cite{cvx}. The P-S algorithm also uses a dual method and
the comparison between it and Algorithm 1 is discussed in Section
\ref{section:dual method}. CVX is a popular convex optimization
solver which uses a core solver SDPT3 or SeDuMi to solve a large
class of convex optimization problems. It is based on the interior
point methods. In the following, we consider three sets of problems.
For each one, we test 30 random instances with input specified in
the following (for problem with size $n=2000$, we only test $3$
instances). Note that the default precision of CVX is $\epsilon =
1.5\times 10^{-8}$. In our dual method, we solve each equation with
precision $\epsilon$. In the gradient projection method, we choose
our starting points to be $0$ and step size to be $1/\sqrt{i}$ in
$i$th iteration, and our stopping criterion is that the objective is
within $\epsilon$ to the CVX optimal value. All the computations are
run on a PC with 1.80GHz CPU and Windows 7 Operating system. We use
CVX Version 1.22 and MATLAB version R2010b. The test results are
shown in Table \ref{table:comparison}.

\begin{table}\centering
\begin{tabular}{||c|c|c||c|c|c|c||}
  \hline\hline
  \# & Problem & $n$ & DM & GP & CVX &
  P-S \\ \hline\hline
  1 & {\bf (TP-1)}  &  50   & 0.050 & 0.214 & 0.378 & 0.154 \\
  2 & {\bf (TP-1)}  & 150   & 0.370 & 0.681 & 1.792 & 3.751 \\
  3 & {\bf (TP-1)}  & 500   & 5.614 & 2.559 & 15.76 & 219.8 \\
  4 & {\bf (TP-1)}  & 2000  & 242.3 & 34.58 & 4953.1 & N/A \\
  \hline\hline
  5 & {\bf (TP-2)}  & 50    & 0.018 & 0.093 & 0.264  & $<$ 0.001 \\
  6 & {\bf (TP-2)}  & 150   & 0.136 & 0.176 & 0.547  & $<$ 0.001 \\
  7 & {\bf (TP-2)}  & 500   & 1.911 & 0.806 & 14.04  & 0.0013 \\
  8 & {\bf (TP-2)}  & 2000  & 58.15 & 3.160 & 4798.9 & 0.0026 \\
  \hline\hline
  9 & {\bf (TP-3)}  & 50    &  0.011  & 0.088 & 1.001   & 0.227    \\
  10 & {\bf (TP-3)} & 150   &  0.015  & 0.124 & 4.950   & 0.966    \\
  11 & {\bf (TP-3)} & 500   &  0.210  & 0.374 & 49.44   & 6.210\\
  12 & {\bf (TP-3)} & 2000  &  0.499  & 1.526 & 2350.2  &  91.66   \\
  \hline\hline
\end{tabular}\caption{Performance Comparisons. DM is the dual method developed in Section \ref{section:dual method}, GP is the gradient projection method developed in Section \ref{section:gradient method} and P-S is the algorithm in \cite{padakandla}. N/A means this method can not return the optimal solution in the corresponding case}\label{table:comparison}
\end{table}

The first problem is
\begin{eqnarray}\label{firsttestproblem}
{\bf (TP-1)}\quad\quad\mbox{minimize} & \sum_{i=1}^n
\left(\frac{1}{4}y_i^4 +
v_iy_i\right) \nonumber\\
\mbox{subject to} & \sum_{i=1}^k y_i \ge \sum_{i=1}^k \alpha_i, &
\forall
k = 1,...,n-1\nonumber\\
& \sum_{i=1}^n y_i = \sum_{i=1}^n \alpha_i \nonumber\\
& y_i\ge 0, & \forall i=1,...,n.
\end{eqnarray}
This problem is considered in the numerical tests in
\cite{padakandla}. We use the same setup, that is, we assume all the
parameters $\alpha_i$ and $v_i$ are drawn from i.i.d. uniform
distributions on $[0,1]$, with $v_i$'s sorted in ascending order. In
order to apply our dual method and gradient projection method, we
first perform a transformation as described in Section
\ref{subsec:equivalent_forms}. Note that the equality constraint in
(\ref{firsttestproblem}) can be replaced by an inequality constraint
since the objective of (\ref{firsttestproblem}) is increasing in
$y$. Then we add an artificial upper bound $\beta = \sum_{i=1}^n
\alpha_i$ to all $y_i$'s and define $z_i = \beta- y_i$. After these
transformations, the problem becomes
\begin{eqnarray*}\label{firsttestproblemtransformed}
\mbox{minimize}_{\vec{z}} & \sum_{i=1}^n f_i(z_i) = \sum_{i=1}^n
(\frac{1}{4}(\beta-z_i)^4 + v_i(\beta
- z_i))\nonumber\\
\mbox{subject to} & \sum_{i=1}^k z_i \le \sum_{i=1}^k (\beta -
\alpha_i), & \forall\mbox{ } k = 1,..,n \nonumber\\
& z_i\ge 0, & \forall\mbox{ }  i=1,..,n,
\end{eqnarray*}
which can be solved by both our dual method and the gradient
projection method.

From Table \ref{table:comparison}, we can see that both our
algorithms outperform the P-S algorithm and CVX for this problem.
The reason that the dual algorithms performs better than the P-S
algorithm is explained in Section \ref{section:dual method}.
Particularly, in this case, $g_i^{-1}(x) = (x - v_i)^{1/3}$ and the
equation (\ref{refequation}) does not have a closed form solution.
In such cases, the complexity of solving (\ref{refequation}) is
essentially similar to the complexity of solving (\ref{iter1case1'})
or (\ref{iter2case1'}). And the computational complexity of our dual
algorithm is less than that of the P-S algorithm by an order of $n$.

Note that in this problem
\begin{eqnarray*}
\bar{z}_i = arg\min_{0\le z\le \beta} f_i(z) = \beta.
\end{eqnarray*}
Therefore $L=n$, that means this is already the worst case scenario
for the dual method. Yet it still performs quite well. Because of
the same reason, the gradient projection method works better than
the dual method in this case. And both of them perform better than
CVX significantly.

The second problem is
\begin{eqnarray}\label{secondtestproblem}
{\bf (TP-2)} \quad\quad\mbox{minimize}_{\vec{y}} & \sum_{i=1}^n \frac{v_i}{1-y_i} \nonumber\\
\mbox{subject to} & \sum_{i=1}^k y_i \ge \sum_{i=1}^k \alpha_i, &
\forall
k = 1,...,n-1\nonumber\\
& \sum_{i=1}^n y_i = \sum_{i=1}^n \alpha_i \nonumber\\
& 0 \le y_i\le 1, & \forall i=1,...,n.
\end{eqnarray}
This problem is also considered in \cite{padakandla}. We again use
the same setup, where $\alpha_i$ and $v_i$ are drawn from i.i.d.
uniform distributions on $[0,1]$, with $v_i$'s sorted in ascending
order. Similar to what we have done for ({\bf TP-1}), we replace the
equality constraint with an inequality constraint and define $z_i =
1-y_i$. An equivalent form of (\ref{secondtestproblem}) is then
obtained as follows:
\begin{eqnarray*}
\mbox{minimize}_{\vec{z}} & \sum_{i=1}^nf_i(z_i) = \sum_{i=1}^n \frac{v_i}{z_i} \nonumber\\
\mbox{subject to} & \sum_{i=1}^k z_i \le k - \sum_{i=1}^k \alpha_i,
& \forall
k = 1,...,n\nonumber\\
& 0 \le z_i\le 1, & \forall i=1,...,n.
\end{eqnarray*}
In this case, $g_i^{-1}(x) = \sqrt{-v_i/x}$. As shown in
\cite{padakandla}, there is a closed form solution to
(\ref{refequation}) in this case. Thus the P-S algorithm can solve
this problem very fast. This is indeed observed in Table
\ref{table:comparison}. Note that the performance of both the dual
and gradient projection methods also improve. This is partly because
it is easier to evaluate the function $g^{-1}$ in this case than in
the first problem. Still, the gradient projection method is faster
than the dual method in this case, because we have $\bar{z}_i = 1$
thus $L=n$ in the dual method. Again, both algorithms work much
faster than CVX.

The last problem is the inventory control problem described in
Section \ref{subsec:application}. The optimization problem is:
\begin{eqnarray*}\label{inventory2}
{\bf (TP-3)} \quad\mbox{minimize}_{\vec{y}} & \sum_{i=1}^n
\left(u_i{\mathbb E}(D_i-y_i)^+ + o_i
{\mathbb E}(y_i-D_i)^+\right)\nonumber\\
\mbox{subject to} & \sum_{i=1}^k y_i \le \sum_{i=1}^k \alpha_i, &
\forall k
= 1,...,n\nonumber\\
& y_i\ge 0, & \forall i = 1,...,n.
\end{eqnarray*}
In the numerical experiments, we assume that $o_i \sim U[5,10]$,
$u_i \sim U[20,25]$ and $\alpha_i \sim U[0,20]$. We also assume that
each $D_i$ follows an exponential distribution with parameter
$\eta_i \sim U[0.1,0.2]$. By applying the property of exponential
distribution, the objective can be equivalently written as:
\begin{eqnarray*}
\sum_{i=1}^n f_i(y_i) = \sum_{i=1}^n
\left(\frac{u_i+o_i}{\eta_i}e^{-\eta_iy_i} + o_iy_i\right)
\end{eqnarray*}
with $\bar{y}_i =
\frac{1}{\eta_i}\log{\left(\frac{u_i+o_i}{o_i}\right)}\in [5.49,
17.92]$, $g_i(y_i) = f'(y_i) = - (u_i+o_i)e^{-\eta_i y_i}+o_i$ and
$g_i^{-1}(s) =
-\frac{1}{\eta_i}\log{\left(\frac{o_i-s}{u_i+o_i}\right)}$. Clearly
(\ref{refequation}) doesn't have a closed form solution with such
$g^{-1}$, therefore our algorithms outperform the P-S algorithm. In
fact, for this problem, the dual method works very well. This is
because the $L$'s in this case are usually much smaller than $n$. In
fact, $L$ is less than $n/4$ in most test problems. This will
greatly reduce the computations in the dual method and make it very
efficient. The gradient method could not take advantage of this
structure and thus only has similar performance as in other
problems.

To summarize the numerical results, we observe that our algorithms
exhibit significant performance improvement over the P-S algorithm
when equation (\ref{refequation}) does not have a closed form
solution. And they also greatly improve over the performance of CVX.
Between the dual method and the gradient projection method, the
former one is more efficient when $L$ is relatively small,
otherwise, the latter one is usually more efficient.

\section{Conclusions}
\label{section:conclusion}

In this paper, we propose two algorithms for solving a class of
convex optimization problems with linear ascending constraints. When
the objective is separable, we propose a dual method which improves
the worst case complexity of the algorithm proposed in
\cite{padakandla}. Furthermore, we propose a gradient projection
algorithm in which each projection step uses the dual method as a
subroutine. The gradient projection algorithm can be used to solve
more general non-separable problems and does not need to evaluate
the inverse of the gradient function which the dual methods usually
require. Numerical results show that both of our proposed algorithms
work well in test problems.

\section{Acknowledgement}
\label{section:acknoledgement} The author thanks Diwakar Gupta and
Shiqian Ma for useful discussions and Arun Padakandla and Rajesh
Sundaresan for sharing the code of the P-S algorithm.

\appendix

\section{Proof of Lemma \ref{lemma1}}
\label{appendix:proofoflemma1}

{\noindent\bf Proof of Lemma \ref{lemma1}.} We prove by
contradiction. If $\tilde{\lambda}_{w_{r^*}} = 0$, then
\begin{enumerate}
\item If there exists $r' = \max\{0 < r < r^*:\tilde{\lambda}_{w_r} > 0
\}$, then we know that
\begin{eqnarray*}
\sum_{s=w_{r'}+ 1}^{w_{r^*}} (\tilde{y}_s-\alpha_s) \le 0.
\end{eqnarray*}
Also, since $r'<r^*$, we know that
\begin{eqnarray*}
\sum_{s=w_{\bar{j}-1}+1}^{w_{r'}}\left(H_s\left(-
\sum_{t=r'+1}^L\tilde{\lambda}_{w_t}\right)-\alpha_s\right) < 0.
\end{eqnarray*}
Therefore, we have
\begin{eqnarray}\label{17}
&&\sum_{s=w_{\bar{j}-1}+1}^{w_{r^*}}\left(H_s\left(-
\sum_{t=r^*+1}^L\tilde{\lambda}_{w_t}\right)-\alpha_s\right)\nonumber\\
&= & \sum_{s=w_{\bar{j}-1}+1}^{w_{r'}}\left(H_s\left(-
\sum_{t=r'+1}^L\tilde{\lambda}_{w_t}\right)-\alpha_s\right) +
\sum_{s=w_{r'}+1}^{w_{r^*}}\left(H_s\left(-
\sum_{t=r^*+1}^L\tilde{\lambda}_{w_t}\right)-\alpha_s\right)\nonumber\\
&= & \sum_{s=w_{\bar{j}-1}+1}^{w_{r'}}\left(H_s\left(-
\sum_{t=r'+1}^L\tilde{\lambda}_{w_t}\right)-\alpha_s\right) +
\sum_{s = {w_{r'}+1}}^{w_{r^*}}(\tilde{y}_s-\alpha_s)<0.
\end{eqnarray}
Here the first equality is because of the assumption that
$\tilde{\lambda}_{w_r} = 0 $ for all $r' < r\le r^*$, and the second
equality is because of the induction assumption. However, (\ref{17})
contradicts with the definition of $r^*$.

\item If all $\tilde{\lambda}_{w_r} = 0$ for $r < r^*$. Then we have
\begin{eqnarray*}
\sum_{s=w_{\bar{j} - 1}+1}^{w_{r^*}} (\tilde{y}_s - \alpha_s )
=\sum_{s=w_{\bar{j} - 1}+1}^{w_{\bar{j}}} (\tilde{y}_s - \alpha_s )
+ \sum_{s=w_{\bar{j}}+1}^{w_{r^*}} (\tilde{y}_s - \alpha_s ) < 0,
\end{eqnarray*}
where in the last inequality, the first term is less than $0$ since
the algorithm enters Case 2 in this loop, and the second term is
less than or equal to $0$ due to the induction assumption. This
contradicts with the definition of $r^*$. Therefore we have proved
Lemma \ref{lemma1}. $\hfill \Box$
\end{enumerate}

\bibliographystyle{plain}
\bibliography{ref}

\end{document}